\documentclass[12pt, reqno]{amsart}
\usepackage{amsmath}
\usepackage{amsthm}
\usepackage{amssymb}
\usepackage{enumerate}
\usepackage{mathrsfs}
\usepackage{graphicx}
\usepackage{graphicx}
\usepackage{amssymb,amsfonts}
\usepackage{amsmath}
\usepackage[colorlinks,linkcolor=blue,anchorcolor=blue,citecolor=blue]{hyperref}
\usepackage[numbers,sort&compress]{natbib}
\usepackage{marvosym}
\usepackage{epstopdf}
\usepackage{caption}
\usepackage{multicol}
\usepackage{multirow}
\usepackage{makebox}
\usepackage{subfigure}
\numberwithin{equation}{section}
\newtheorem{theo}{Theorem} 

\newtheorem{lem}{Lemma}

\begin{document}

\title {$\textrm{SL}(n)$ covariant matrix-valued valuations on $L^{p}$-spaces }

\author[Zeng and Lan]{Chunna Zeng$^{*}$ and Yu Lan}

\address{1.School of Mathematical Sciences, Chongqing Normal University, Chongqing 401331, People's Republic of China; Institut f\"{u}r Diskrete Mathematik und Geometrie, Technische Universit\"{a}t Wien, Wiedner Hauptstra\ss e 8--10/1046, 1040 Wien, Austria }
\email{zengchn@163.com}

\address{2.School of Mathematics and Sciences, Chongqing Normal University, Chongqing 401332, People's Republic of China}
\email{18181567221@163.com}

\thanks{The first author is supported by the Major Special Project of NSFC (Grant No. 12141101), the Young Top-Talent program of Chongqing (Grant No. CQYC2021059145), the Science and Technology Research Program of Chongqing Municipal Education Commission (Grant No. KJZD-K202200509) and Natural Science Foundation Project of Chongqing (Grant No. CSTB2024NSCQ-MSX0937)}
\thanks{{*}Corresponding author: Chunna Zeng}
\begin{abstract}
A complete classification is established for continuous and $\textrm{SL}(n)$ covariant matrix-valued valuations on $L^{p}(\mathbb{R}^n,\vert x\vert^2dx)$. The assumption of matrix symmetry  is eliminated. For $n\geq 3$, such valuation is uniquely characterized by the moment matrix of  measurable function. In the 2-dimensional case, while the rotation matrix shows up.
\end{abstract}

\maketitle

\section{Introduction}
According to  Erlangen program proposed by Felix Klein in 1872, the research and  classification of geometric notions compatible with
transformation groups are significant issues in geometry. Numerous operators defined on geometric entities adhere to the inclusion-exclusion principle, making it natural to consider the property of being a valuation in classification. A function $\mathrm{Y}: \mathcal{W}\rightarrow \mathcal{G}$ mapping a collection
$\mathcal{W}$ of sets to an abelian semigroup $\langle \mathcal{G}, +\rangle$
is  termed a valuation if it satisfies
\begin{equation}
\mathrm{Y}(P_{1})+\mathrm{Y}(P_{2})=\mathrm{Y}(P_{1}\cup P_{2})+\mathrm{Y}(P_{1}\cap P_{2}),
\end{equation}
for all $P_{1}, P_{2}$ with $ P_{1}\cup P_{2}, P_{1}\cap P_{2}\in\mathcal{W}.$ The idea of valuations on convex bodies emerged from Dehn's solution of Hilbert's third problem.
The Hadwiger characterization theorem is a landmark result in this field, providing a comprehensive classification of continuous, rigid motion invariant valuations on convex bodies.
This theorem simplifies these proofs for many results in geometric probability and integral geometry theory.
And it is the cornerstone for many beautiful results in modern valuation theory. Professor Ludwig has made significant contributions to  body-valued valuations \cite{Ludwig1,Ludwig2,Ludwig3,Ludwig4,Ludwig9}. For  historical development and recent contributions on convex bodies, star shaped sets, manifolds and so on, see \cite{Alesker1,Alesker2,Alesker3,Alesker4,Klain1,Klain2,Klain3,Klain-Rota,L-Y-Z1,Ludwig5,Ludwig6,Ludwig7,Ludwig8,McMullen1,McMullen2,Ma,Parapatits1,Parapatits2,Wannerer,Wang1,Zeng-Ma}.

\par Recent research has increasingly focused on valuations in function spaces.
Assume that $\mathcal{S}$ is  a lattice of some functions. For $h, f\in\mathcal{S},$
the operations $h\vee f$ and $h\wedge f$ denote $\max\{h, f\}$ and $\min\{h, f\}$, respectively.
A function $\Psi$ mapping  $(\mathcal{S},\vee,\wedge)$ to an abelian semigroup  is referred to as a valuation if it satisfies
\begin{equation}\label{0}
\Psi(h\vee f)+\Psi(h\wedge f)=\Psi(h)+\Psi(f)
\end{equation}
for all $h, f, h\vee f, h\wedge f\in \mathcal{S}.$
Valuations are extensively studied in $L^{p}$-spaces \cite{Ludwig6,Tsang1,Tsang2,Wang2}, Sobolev spaces \cite{Ludwig5,Ludwig7}, Lipschitz functions \cite{Colesanti1,Colesanti2}, and spaces of functions of bounded variation \cite{Wang}, among others.

\par In particular, the theory of valuations in $L^{p}$-spaces has particularly attracted the attention of geometers.
For instance, Tsang \cite{Tsang1} first classified continuous, translation invariant real-valued valuation on $L^{p}(\mathbb R^{n})$, and also classified  continuous, rotation invariant valuation on $L^{p}(S^{n-1})$.
Building on Tsang's work, Wang \cite{Wang3} extended these results to vector-valued valuations,  showing that all continuous $\textrm {SL}(n)$ covariant vector-valued valuations on $L^{p}(\mathbb{R}^n, \vert x\vert dx)$ are  moment vectors of measurable functions. In 2013, Ludwig \cite{Ludwig6} classified  continuous $\textrm{SL}(n)$ covariant symmetric matrix-valued valuations on $L(\mathbb{R}^n,\vert x\vert^{2} dx)$,  demonstrating that such valuations are  moment matrices. Wang \cite{Wang2} later extended this classification to  continuous, $\textrm{SL}(n)$ covariant  matrix-valued valuations on $L(\mathbb{R}^n,\vert x\vert^{2} dx)$  without the assumption of matrix symmetry.

\par
\par For $p\in[1,\infty)$, let $L^{p}(\mathbb{R}^n, \vert x\vert^{2}dx)$ denote the space of measurable functions $h:\mathbb{R}^n\rightarrow\mathbb{R}$ satisfying that $\int_{\mathbb{R}^n}\vert h(x)\vert^{p}\vert x\vert^{2}dx<\infty$. Here, $\vert x\vert$ represents the norm of $x$ in $\mathbb{R}^n$.
Denote by $\mathbb{M}^n$  the set of $n\times n$ matrix in $\mathbb{R}^n$, and
 $\mathbb{M}^n_{e}$  its subset consisting of  $n\times n$  symmetric matrix.
An operator $\Psi: L^{p}(\mathbb{R}^n,\vert x\vert^2dx)\rightarrow\mathbb{M}^n$ is termed $\textrm{SL}(n)$ covariant if for all $h\in L^{p}(\mathbb{R}^n,\vert x\vert^2dx)$ and $\phi\in\textrm{SL}(n)$,
\begin{equation*}
\Psi(h\circ\phi^{-1})=\phi \Psi(h)\phi^{t}.
\end{equation*}
The moment matrix $\mathrm{K}(h)$ for  a measurable function $h:\mathbb{R}^n\rightarrow\mathbb{R}$
  is the $n\times n$ matrix with  entries (which may be infinite),
\begin{equation*}
\mathrm{K}_{ij}(h)=\int_{\mathbb{R}^n} h(x) x_{i}x_{j}dx.
\end{equation*}

Denote by $C(\mathbb{R})$ the collection of continuous functions $\xi:\mathbb{R}\rightarrow\mathbb{R}.$  Denote by  $A(\mathbb{R})$ the collection of continuous functions $\xi:\mathbb{R}\rightarrow\mathbb{R}$ satisfying that there exits $d\in\mathbb{\mathbb{R}}$ for which $\vert\xi(t)\vert\leq d\vert t\vert\ $
 holds for any $t\in\mathbb{R}$.
Similarly, $\widetilde{A}(\mathbb{R})$
 is the collection of continuous functions
$\xi:\mathbb{R}\rightarrow\mathbb{R}$ satisfying that there exists   $d\in\mathbb{\mathbb{R}}$ for which $\vert\xi(t)\vert\leq d{\vert t\vert}^{p} $ for every $t\in\mathbb{R}$.
Professor Lugwid provided  the following theorem.

\begin{theo}{\cite{Ludwig6}} Suppose that $n\geq3$.
An operator 
$\Psi: L(\mathbb{R}^n,\vert x\vert^{2} dx)\rightarrow\mathbb{M}^n_{e}$
is a continuous $\mathrm{SL}(n)$ covariant valuation if and only if there exists $\xi\in A(\mathbb{R})$ such that
\begin{equation*}
\Psi(h)=\mathrm{K}(\xi\circ h)
\end{equation*}
for all $h\in L(\mathbb{R}^n,\vert x\vert^{2} dx)$.
\end{theo}

Recently, Wang \cite{Wang2} extended Ludwig's result \cite{Ludwig6} by removing the symmetry assumption. Additionally, he provided
the classification of the 2-dimensional case.

\begin{theo}\label{theo1}\cite{Wang2}
An operator 
$\Psi:L(\mathbb{R}^2,\vert x\vert^{2}dx)\rightarrow\mathbb{M}^2$
is a continuous $\mathrm{SL}(2)$ covariant valuation if and only if there exist $\xi\in A(\mathbb{R})$ and a constant $s\in\mathbb{R}$ such that
\begin{equation*}
\Psi(h)=\mathrm{K}(\xi\circ h)+s\rho_{\frac{\pi}{2}}
\end{equation*}
for every $h\in  L(\mathbb{R}^2,\vert x\vert^{2} dx)$. Here $\rho_{\frac{\pi}{2}}$ represents the counter-clockwise rotation  by $\pi/2$ in $\mathbb{R}^2$.
\end{theo}

\begin{theo}\label{theo2}\cite{Wang2}
Suppose that $n\geq3$. An operator 
$\Psi: L(\mathbb{R}^n,\vert x\vert^{2} dx)\rightarrow\mathbb{M}^n$ is a continuous $\mathrm{SL}(n)$ covariant valuation if and only if there exists  $\xi\in A(\mathbb{R})$ such that
\begin{equation*}
\Psi(h)=\mathrm{K}(\xi\circ h)
\end{equation*}
for all $h\in L(\mathbb{R}^n,\vert x\vert^{2} dx)$.
\end{theo}

The purpose of our paper is to characterize continuous  and $\textrm{SL}(n)$ covariant matrix-valued valuations on $L^{p}(\mathbb{R}^n, \vert x\vert^{2}dx)$. We will extend Wang's result (\cite{Wang2}) from the space $L(\mathbb{R}^n,\vert x\vert^2dx)$  to $L^{p}(\mathbb{R}^n,\vert x\vert^2dx)$.

\begin{theo}\label{theo3}
An operator $\Psi:L^{p}(\mathbb{R}^2,\vert x\vert^{2}dx)\rightarrow\mathbb{M}^2$ is a continuous  $\mathrm{SL}(2)$ covariant valuation if and only if there exist $\xi\in \widetilde{A}(\mathbb{R})$ and a constant $s\in\mathbb{R}$ such that
\begin{equation*}
\Psi(h)=\mathrm{K}(\xi\circ h)+s\rho_{\frac{\pi}{2}}
\end{equation*}
for every $h\in L^{p}(\mathbb{R}^2,\vert x\vert^{2}dx)$. Here $\rho_{\frac{\pi}{2}}$ represents the counter-clockwise rotation  by $\pi/2$ in $\mathbb{R}^2$.
\end{theo}

\begin{theo}\label{theo4}
Suppose that $n\geq3$. A function $\Psi: L^{p}(\mathbb{R}^n,\vert x\vert^{2}dx)\rightarrow\mathbb{M}^n$ is a continuous $\mathrm{SL}(n)$ covariant valuation if and only if there exists  $\xi\in \widetilde{A}(\mathbb{R})$ such that
\begin{equation*}
\Psi(h)=\mathrm{K}(\xi\circ h)
\end{equation*}
for all $h\in L^{p}(\mathbb{R}^n,\vert x\vert^{2}dx)$.
\end{theo}

\section{~Preliminaries}
\par Consider an arbitrary set $X$  and  a $\sigma$-algebra $\mathcal{A} $ on $X$. Assume $\mu: \mathcal{A}\rightarrow [0, \infty]$ is a measure defined on $\mathcal{A} $. Then $(X,\mathcal{A},\mu)$ forms a measurable space. The space $L^{p}(\mu)$, for $ p\in[1,\infty),$ consists of $\mu$-measurable functions $h: X\rightarrow[-\infty, +\infty]$ satisfying that
\begin{equation*}
\int_{X}\lvert h \rvert^{p}d\mu<\infty.
\end{equation*}
Denote by $\Vert h\Vert_{p}$ the $L^{p}(\mu)$-norm of $h,$ where $h\in L^{p}(\mu).$ The
 expression of $\Vert h\Vert_{p}$ is

\begin{equation*}
\Vert h\Vert_{p}=\left(\int_{X}\lvert h \rvert^{p}d\mu\right)^{\frac{1}{p}}.
\end{equation*}
In this case, the functional $\Vert\cdot\Vert_{p}: L^{p}(\mu)\rightarrow\mathbb{R}$ acts as a semi-norm. Since the functions in $L^{p}(\mu)$ that are equal a.e. with respect to $\mu$ are considered identical, which converts
$\Vert\cdot\Vert_{p}$ into a norm, thus making $L^{p}(\mu)$ a normed linear space.


\par In this paper, we investigate two measurable spaces: $(\mathbb{R}^n,\mathscr{B},\lambda)$ and $(\mathbb{R}^n,\mathscr{B},\mu_{n})$. Here,  $\mathscr{B}$ represents the set of Lebesgue measurable sets in $\mathbb{R}^n$, $\lambda$ denotes the Lebesgue measure, and $\mu_{n}$ is  defined by
\begin{equation*}
\mu_{n}(D)=\int_{D}\vert x \vert^{2} dx
\end{equation*}
for $D\in\mathscr{B}$. We write $L^{p}(\mu_{n})$ as $L^{p}(\mathbb{R}^n,\vert x\vert^{2}dx)$. And we use $L^{p}(\mathbb{R}^n)$ as a replacement for
 $L^{p}(\lambda)$, meaning that for
  $h\in L^{p}(\mathbb{R}^n)$, we typically write $\int_{\mathbb{R}^n}h(x)dx$  instead of $\int_{\mathbb{R}^n}hd\lambda$. Likewise, for $h\in L^{p}(\mu_{n})$, we replace $\int_{\mathbb{R}^n}h d\mu_{n} $ by $\int_{\mathbb{R}^n}h(x)\vert x \vert^2 dx$. For convenience, we denote either $\lambda$-measurable functions or $\mu_{n}$-measurable as the measurable functions. Notice that if $D\in\mathscr{B}$, then $\mu_{n}(D)=0$ and $\lambda(D)=0$ are equivalent, so we abbreviate almost everywhere with respect to $\mu_{n}$ and  $\lambda$ as a.e.. The $L^{p}$-norm of $h$ in $L^{p}(\mu_{n})$ is expressed as
\begin{equation*}
\Vert h\Vert_{L^{p}(\mu_{n})}=\bigg(\int_{\mathbb{R}^n}\lvert h(x) \rvert^{p}\vert x \vert^{2}dx\bigg)^{\frac{1}{p}}.
\end{equation*}
If $\Vert h_{k}-h\Vert_{L^{p}(\mu_{n})}\rightarrow 0$ as $k\rightarrow\infty$, then we call $h_{k}\rightarrow h$ as $k\rightarrow\infty$ in $L^{p}(\mu_{n})$. A function $\Psi:L^{p}(\mu_{n})\rightarrow\mathbb{M}^n$ is called continuous if $\Psi(h_{k})\rightarrow \Psi(h)$ , as $h_{k}\rightarrow h$ in $L^{p}(\mu_{n})$.
On the other hand, we remain use  $\Vert h\Vert_{P}$ to denote the $L^{p}$-norm of $h$ in  $ L^{p}(\mathbb{R}^n)$,
\begin{equation*}
\Vert h\Vert_{p}=\bigg(\int_{\mathbb{R}^n}\vert h\vert^pdx\bigg)^\frac{1}{p}.
\end{equation*}
Denote by $1_{P}$ the indicator function of $P\subset\mathbb{R}^n$. It shows that ${1}_{P}(x)=1$ if $x\in P$ and $1_{P}(x)=0 $ otherwise.

\par Assume that $e_{1}, \cdots, e_{n}$ is the standard basis  in Euclidean space $\mathbb{R}^n.$
For $\textbf{a}=(a_{1},\cdots,a_{n})\in\mathbb{R}^{n}$ and $\textbf{b}=(b_1,\cdots, b_{n})\in\mathbb{R}^{n}$, the scalar product  is given by  $\textbf{a}\cdot \textbf{b}=a_{1}b_{1}+\cdots+a_{n}b_{n}$.
Let $\mathcal{P}^n$ represent the set of compact convex polytopes, and $\mathcal{P}^{n}_{0}$ its subset consisting of
  compact convex polytopes that contains the origin $0$. Recall that both spaces are endowed with the topology induced by the Hausdorff metric.
 For $P\in\mathcal{P}^{n} $, the moment matrix $\mathrm{M}(P)$ is the $n\times n$ matrix with coefficients given by
\begin{equation*}
\mathrm{M}_{ij}(P)=\int_{P}x_{i}x_{j}dx.
\end{equation*}
A function $\mathrm{Y}:\mathcal{P}^n\rightarrow \mathbb{M}^{n}$ is said to be  $\mathrm{SL}(n)$ covariant if it satisfies
\begin{equation*}
\mathrm{Y}(\phi P)=\phi \mathrm{Y}(P)\phi^{t}
\end{equation*}
for any $P\in \mathcal{P}^n$ and $\phi\in \mathrm{SL}(n)$. In addition, the moment matrix  is $\mathrm{SL}(n)$ covariant.

\par A function $\mathrm{Y}:\mathcal{P}\rightarrow\mathbb{M}^n$ is termed weakly simple if there exists a matrix $M_{0}$ such that $\mathrm{Y}(P)=M_{0}$ for any $P\in\mathcal{P}^n$ where  dim$ P$ is at most $n-1$ (see \cite{{Wang2}}). If $M_{0}=\mathbf{0}$, then we say the function $\mathrm{Y}$ is simple.
The following quotes are the conclusion of the real analysis.

\begin{lem}\label{lem1}\cite{Rudin}
The function $\Phi(h)=\vert h\vert^{p}$, mapping $L^{p}(\mu_{n})$ to $ L^{1}(\mu_{n})$
 for $h\in L^{p}(\mu_{n})$, is continuous.
\end{lem}

\par A closed cube $C$ in $\mathbb{R}^n$ is defined as
\begin{equation*}
C=\{(x_{1},\cdots,x_{n}): c_{i}\leq x_{i}\leq d_{i},\  x_{i}, c_{i}, d_{i}\in\mathbb{R}\},
\end{equation*}
where $\vert d_{i}-c_{i}\vert= \vert c_{j}-d_{j}\vert$,  $1\leq i, j\leq n$. For subsets $P, Q\subseteq\mathbb{R}^n$, the definition of the symmetric difference of  $P, Q$ is
 \begin{equation}
 P\bigtriangleup Q=P\setminus Q \cup Q\setminus P.
 \end{equation}

\begin{lem}\label{lem2}\cite{Tsang1}
Suppose that $\varepsilon>0$ and $D\in\mathscr{B}$ with $\lambda(D)<\infty$, then there exists a series of closed cubes $C_{1}, C_{2},\ldots, C_{k}\subseteq\mathbb{R}^n$ such that $\lambda(D\bigtriangleup \bigcup\limits_{i=1}\limits^{k} C_{i})<\varepsilon$.
\end{lem}
\begin{lem}\label{lem3}
Let $\Psi:L^{p}(\mu_{n})\rightarrow\mathbb{M}^2$ be an $\mathrm{SL}(2)$ covariant function. Then there exists a constant $s\in\mathbb{R}$ such that
\begin{equation*}
\Psi(0)=s\rho_{\frac{\pi}{2}}.
\end{equation*}
 Here $\rho_{\frac{\pi}{2}}$ represents the counter-clockwise rotation  by $\pi/2$ in $\mathbb{R}^2$.
\end{lem}

\par \emph{Proof.} \
Define two linear transforms $\phi_{1}, \phi_{2}\in \textrm{SL}(2)$ by  $\phi_{1}$=
$\begin{pmatrix}
k & 0\\
0 & k^{-1}
\end{pmatrix}$$(k\neq0),$
and $\phi_{2}=$
$\begin{pmatrix}
1 & 0\\
1 & 1
\end{pmatrix}$.
The  $\textrm{SL}(2)$ covariance of $\Psi$ shows that
\begin{equation}\label{6}
 \Psi(0)= \Psi(0\circ\phi_{i}^{-1})=\phi_{i}\Psi(0)\phi_{i}^{t},
\end{equation}
 where $i=1,2$.
Assume that $\Psi(0)$=
$\begin{pmatrix}
a_{11} & a_{12}\\
a_{21} & a_{22}
\end{pmatrix}$,
 taking $\phi_1$ in (\ref{6}) we obtain that
 \begin{equation*}
\begin{pmatrix}
 a_{11} & a_{12}\\
a_{21} & a_{22}
 \end{pmatrix}=
 \begin{pmatrix}
 k^{2}a_{11} & a_{12}\\
 a_{21} & k^{-2}a_{22}
 \end{pmatrix}
 \end{equation*}
 for $k\neq0$. Then we have $a_{11}=0$ and $a_{22}=0$.
 Taking $\phi_2$ in (\ref{6}) we have that
\begin{equation*}
\begin{pmatrix}
 0 & a_{12}\\
a_{21} & 0
 \end{pmatrix}=
 \begin{pmatrix}
 0 & a_{12}\\
 a_{21} & a_{12}+a_{22} \end{pmatrix},
 \end{equation*}
 thus $a_{12}+a_{21}=0.$
 \qed
\begin{lem}\label{lem4}
Suppose that $n\geq3$. A operator $\Psi:L^{p}(\mu_{n})\rightarrow\mathbb{M}^n$ is  an $\mathrm{SL}(n)$ covariant, then
\begin{equation*}
\Psi(0)=\mathrm{\textbf{0}}.
\end{equation*}
\end{lem}

\emph{Proof.} \
For a real number $k\neq 0$, let
\begin{equation*}
\varphi_{n}=
\begin{pmatrix}
k &   & &\\
  & k &  &\\
  &   & \ddots &\\
  &   &  & k^{1-n}
\end{pmatrix}.
\end{equation*}
It is obvious that $\varphi_{n}\in\mathrm{SL}(n)$. Since that $\Psi$  has the property of  $\mathrm{SL}(n)$ covariant, then
\begin{equation}\label{7}
 \Psi(0)= \Psi(0\circ\varphi_{n}^{-1})=\varphi_{n}\Psi(0)\varphi_{n}^{t}.
\end{equation}
Set
\begin{equation*}
\Psi(0)=
\begin{pmatrix}
a_{11} & a_{12} & \cdots & a_{1n}\\
a_{21} & a_{22} & \cdots & a_{1n}\\
\vdots & \vdots & \ddots & \vdots\\
a_{n1} & a_{n2} & \cdots & a_{nn}
\end{pmatrix},
\end{equation*}
thus (\ref{7}) becomes
\begin{equation*}
\begin{pmatrix}
a_{11} & a_{12} & \cdots & a_{1n}\\
a_{21} & a_{22} & \cdots & a_{1n}\\
\vdots & \vdots & \ddots & \vdots\\
a_{n1} & a_{n2} & \cdots & a_{nn}
\end{pmatrix}=
\begin{pmatrix}
k^{2}a_{11} & k^{2}a_{12} & \cdots & k^{2}a_{1,n-1} & k^{2-n}a_{1n}\\
k^{2}a_{21} & k^{2}a_{22} & \cdots & k^{2}a_{2,n-1} & k^{2-n}a_{2n}\\
\vdots & \vdots & \ddots & \vdots & \vdots\\
k^{2}a_{n-1,1} & k^{2}a_{n-1,2} & \cdots & k^{2}a_{n-1,n-1} & k^{2-n}a_{n-1,n}\\
k^{2-n}a_{n1} & k^{2-n}a_{n2} & \cdots & k^{2-n}a_{n,n-1} & k^{2-2n}a_{nn}
\end{pmatrix}.
\end{equation*}
 Then
 \begin{equation}
 a_{ij}=k^{2}a_{ij},\ \ a_{in}=k^{2-n}a_{in}, \ \  a_{ni}=k^{2-n}a_{ni},\ \ \ 1 \leq i, j\leq n-1;
\end{equation}
and
 \begin{equation}
a_{nn}=k^{2-2n}a_{nn}.
\end{equation}
Therefore,
 \begin{equation}
  a_{ij}=0,\ \ a_{in}=a_{ni}=0,\ \ a_{nn}=0,\ \ 1 \leq i, j\leq n-1.
\end{equation}
It leads to
$\Psi(0)=\textbf{0}.$
\qed
\section{Proof of main theorems}
\begin{lem}\label{lem7} Suppose that $n\geq2$ and  $\xi\in\widetilde{A}(\mathbb{R})$, then the operator $\Psi:L^{p}(\mu_{n})\rightarrow\langle\mathbb{M}^n,+\rangle$ defined by
\begin{equation*}
\Psi(h)=\mathrm{K}(\xi\circ h),
\end{equation*}
is a continuous $\mathrm{SL}(n)$ covariant valuation.
\end{lem}
\emph{Proof.} Because that  $\vert \xi(t)\vert\leq d\vert t\vert^{p}$ for any $t\in\mathbb{R}$, thereby
\begin{align*}
\vert \mathrm{K}(\xi\circ h)\vert&=\vert\int_{\mathbb{R}^{n}}xx^{t}\vert \xi(h(x))\vert dx\vert\\
&\leq d\int_{\mathbb{R}^n} \vert h(x)\vert^{p}\vert x\vert^{2} dx.
\end{align*}
 In addition, since $h\in L^{p}(\mu_{n})$, it can be concluded that $\mathrm{K}(\xi\circ h)<\infty$.

\par Now, we intend to demonstrate that $\Psi$ is a valuation in $L^{p}(\mu_{n})$. Let $h,f\in L^{p}(\mu_{n})$, then
\begin{align*}
\Psi(h\vee f)+\Psi(h\wedge f)&=\int_{\mathbb{R}^n}xx^{t}\vert\xi\circ(h\vee f)\vert dx+\int_{\mathbb{R}^n}xx^{t}\vert\xi\circ(h\wedge f)\vert dx\\
&=\int_{\{h\geq f\}}xx^{t}\vert\xi\circ h\vert dx+\int_{\{h<f\}}xx^{t}\vert\xi\circ f\vert dx\\
&+\int_{\{h\geq f\}}xx^{t}\vert\xi\circ f\vert dx+\int_{\{h<f\}}xx^{t}\vert\xi\circ h\vert dx\\
&=\int_{\mathbb{R}^n}xx^{t}\vert\xi\circ h\vert dx+\int_{\mathbb{R}^n}xx^{t}\vert\xi\circ f\vert dx\\
&=\Psi(h)+\Psi(f).
\end{align*}
Hence $\Psi$ is a valuation on $L^{p}(\mu_{n})$.
\par Next we prove that $\Psi$ is continuous. Suppose that $h_{k}\rightarrow h$ as $k\rightarrow\infty$ in $L^{p}(\mu_{n})$, then we have
\begin{align*}
\vert \Psi(h_{k})-\Psi(h)\vert&=\vert \mathrm{K}(\xi\circ h_{k})- \mathrm{K}(\xi\circ h)\vert\\
&=\vert\int_{\mathbb{R}^n}xx^{t}\vert\xi\circ (h_{k}(x))\vert dx-\int_{\mathbb{R}^n}xx^{t}\vert\xi\circ( h(x))\vert dx\vert\\
&\leq d\int_{\mathbb{R}^n}\vert\lvert h_{k}(x)\vert^{p}-\vert h(x)\vert^{p}\rvert \lvert x\rvert^{2}dx.
\end{align*}
By Lemma \ref{lem1}  we obtain $\Psi(h_{k})\rightarrow \Psi(h)$ as $h_{k}\rightarrow h$ in $L^{p}(\mu_{n})$. Therefore  $h\mapsto \mathrm{K}(\xi\circ h)$ is continuous in $L^{p}(\mu_{n})$.

To demonstrate that $\Psi$ is $\mathrm{SL}(n)$ covariant, observe that due to the $\mathrm{SL}(n)$ covariance of $\mathrm{K}$, we have
\begin{equation*}
\Psi(h\circ\phi^{-1})=\mathrm{K}(\xi\circ(h\circ\phi^{-1}))=\mathrm{K}((\xi\circ h)\circ\phi^{-1})=\phi \mathrm{K}(\xi\circ h)\phi^{t}=\phi \Psi(h)\phi^{t}
\end{equation*}
for all $\phi\in \mathrm{SL}(n)$. Thus, $h\mapsto \mathrm{K}(\xi\circ h)$ is $\mathrm{SL}(n)$ covariant.
\qed
\begin{lem}\label{lem8}\cite{Wang2} An operator $\mathrm{Y}:\mathcal{P}^2\rightarrow\langle\mathbb{M}^2,+\rangle$ is a continuous, weakly simple and $\mathrm{SL}(2)$ covariant valuation if and only if there exist constants $s_{1}, s_{2}\in\mathbb{R}$ such that
\begin{equation*}
\mathrm{Y}(P)=s_{1}\mathrm{M}(P)+s_{2}\rho_{\frac{\pi}{2}}
\end{equation*}
for all $P\in \mathcal{P}^2$.
\end{lem}

\begin{lem}\label{lem9}\cite{Wang2}
Suppose that  $n\geq3$. An operator $\mathrm{Y}:\mathcal{P}^n\rightarrow\langle\mathbb{M}^n,+\rangle$ is a continuous, simple and $\mathrm{SL}(n)$ covariant valuation if and only if there exists a constant $s\in\mathbb{R}$ such that
\begin{equation*}
\mathrm{Y}(P)=s\mathrm{M}(P)
\end{equation*}
for all $P\in\mathcal{P}^n$.
\end{lem}

The  lemma below establishes a significant connection between  functions in $L^{p}(\mu_{n})$ and functions on $\mathcal{P}^n$.
\begin{lem}\label{lem10}\cite{Ludwig6}
For $P\in\mathcal{P}^n$ and $\alpha\in\mathbb{R}$, we have $\mathrm{K}(\alpha1_{P})=\alpha \mathrm{M}(P)$.
\end{lem}

\begin{lem}\label{lem11}
If the operator $\Psi:L^{p}(\mu_{n})\rightarrow\langle\mathbb{M}^2,+\rangle$ is a continuous, $\mathrm{SL}(2)$ covariant valuation, then there exist a constant $s\in\mathbb{R}$ and  $\xi\in C(\mathbb R)$ such that
\begin{equation*}
\Psi(\alpha1_{P})=\xi(\alpha)\mathrm{K}(\alpha1_{P})+s\rho_{\frac{\pi}{2}}
\end{equation*}
for all $P\in\mathcal{P}^2$.
\end{lem}
\emph{Proof.} For some $\alpha\in\mathbb{R}$, define $\mathrm{Y}_{\alpha}:\mathcal{P}^2\rightarrow\langle\mathbb{M}^2,+\rangle$ by setting
\begin{equation*}
\mathrm{Y}_{\alpha}(P)=\Psi(\alpha1_{P}).
\end{equation*}
Because that $\Psi$ is a valuation in $L^{p}(\mu_{n})$, it follows for $P_{1},P_{2},P_{1}\cup P_{2},P_{1}\cap P_{2}\in\mathcal{P}^2$ that
\begin{align*}
\mathrm{Y}_{\alpha}(P_{1}\cup P_{2})+\mathrm{Y}_{\alpha}(P_{1}\cap P_{2}) &=\Psi(\alpha1_{P_{1}\cup P_{2}})+\Psi(\alpha1_{P_{1}\cap P_{2}})\\
&=\Psi(\alpha1_{P_{1}}\vee \alpha1_{P_{2}})+\Psi(\alpha1_{P_{1}}\wedge \alpha1_{P_{2}})\\
&=\mathrm{Y}_{\alpha}(P_{1})+\mathrm{Y}_{\alpha}(P_{2}).
\end{align*}
Therefore, $\mathrm{Y}_{\alpha}:\mathcal{P}^2\rightarrow \langle\mathbb{M}^2,+\rangle$ is a valuation.

 According to Lemma \ref{lem2}, for $P\in\mathcal{P}^2$, there is a series of  $P_{i}$ be a  union of finite closed cubes satisfying $\lambda(P_{i}\bigtriangleup P)<{1}/{i}$. For $\alpha\in\mathbb{R}$, if $x\in P_{i}\cup P$, then there exists a constant $a>0$ such that $\vert x\vert\leq a.$ Notice that
\begin{align*}
\Vert\alpha1_{P_{i}}-\alpha1_{P}\Vert_{L^{p}(\mu_{n})}
 &=\vert\alpha\vert
\bigg(\int_{\mathbb{R}^2}\vert1_{P_{i}}(x)-1_{P}(x)\vert^{p}\vert x\vert^{2} dx\bigg)^{\frac{1}{p}}\\
&\leq\vert\alpha\vert\bigg(\int_{\mathbb{R}^2}\vert1_{P_{i}}(x)-1_{P}(x)
\vert^{p}a^{2}dx\bigg)^{\frac{1}{p}}\\
 &=\vert \alpha\vert a^{\frac{2}{p}}\lambda(P_{i}\bigtriangleup P)^{\frac{1}{p}}.
\end{align*}
So if $\lambda(P_{i}\bigtriangleup P)\rightarrow 0$ as $i\rightarrow\infty,$ then  $\Vert\alpha1_{P_{i}}-\alpha1_{P}\Vert_{L^{p}(\mu_{n})}\rightarrow0$. The continuity of $\Psi$ gives $\Psi(\alpha1_{P_{i}})\rightarrow \Psi(\alpha1_{P})$. Hence, $\mathrm{Y}_{\alpha}(P_{i})\rightarrow \mathrm{Y}_{\alpha}(P)$. It leads that $\mathrm{Y}_{\alpha} $ is continuous.
\par Next we demonstrate that $\mathrm{Y_{\alpha}}$ is $\mathrm{SL}(2)$ covariant. Because of  the $\mathrm{SL}(2)$ covariance of $\Psi$, then for every $\phi\in \mathrm{SL}(2)$,
\begin{equation*}
\mathrm{Y}_{\alpha}(\phi P)=\Psi(\alpha1_{\phi P})=\Psi(\alpha1_{P}\circ\phi^{-1})=\phi \Psi(\alpha1_{P})\phi^{t}
\end{equation*}
for $P\in\mathcal{P}^2$. It follows that $\mathrm{Y_{\alpha}}$ is $\mathrm{SL}(2)$ covariant.
\par By Lemma \ref{lem3},  there exists a constant $s\in\mathbb{R}$ such that
$\Psi(0)=s\rho_{\frac{\pi}{2}}$. Since $\Vert\alpha1_{P}\Vert_{L^{p}(\mu_{n})}=0$ for any $P\in\mathcal{P}^2$ with dim $(P)<2$, then we have $\alpha1_{P}=0$ a.e. in $L^{p}(\mu_{n})$. The continuity of $\Psi$ implies that
\begin{equation}\label{11}
\Psi(\alpha1_{P})=\Psi(0)=s\rho_{\frac{\pi}{2}}
\end{equation}
for any $P\in\mathcal{P}^2$ with dim $(P)<2$. From the definition of $\mathrm{Y}_{\alpha}$,  we obtain that
\begin{equation*}
\mathrm{Y}_{\alpha}(P)=s\rho_{\frac{\pi}{2}},
\end{equation*}
for any $P\in\mathcal{P}^2$ with dim $(P)<2$.
So $\mathrm{Y}_{\alpha}$ is a weakly simple valuation.

\par Now we  have proved that $\mathrm{Y}_{\alpha}$ is a continuous, weakly simple and $\mathrm{SL}(2)$ covariant valuation. From (\ref{11}), Lemma \ref{lem8} and Lemma \ref{lem10},
there exist a constant $\alpha\in\mathbb{R}$ and a functional $s_{\alpha}:\mathbb{R}\rightarrow\mathbb{R}$ such that
\begin{equation*}
\Psi(\alpha1_{P})=s_{\alpha}\mathrm{K}(1_{P})+s\rho_{\frac{\pi}{2}}
\end{equation*}
for all $P\in\mathcal{P}^2$. Define the function $\xi(\alpha)=s_{\alpha}$
for $\alpha\in\mathbb{R}$.
Next, we demonstrate the continuity of $\xi$. Consider $\alpha\in\mathbb{R}$ and a sequence $\{\alpha_{i}\}$ $\subseteq \mathbb{R}$ such that $\alpha_{i}\rightarrow\alpha$. We
observe the following
\begin{align*}
\Vert\alpha_{i}1_{P}-\alpha1_{P}\Vert_{L^{p}(\mu_{n})}=\vert\alpha_{i}
-\alpha\vert\bigg(\int_{\mathbb{R}^2}\vert1_{P}(x)\vert^{p}\vert x \vert^{2} dx\bigg)^{\frac{1}{p}}.
\end{align*}
Thus $\Vert\alpha_{i}1_{P}-\alpha1_{P}\Vert_{L^{p}(\mu_{n})}\rightarrow0$ as $\alpha_{i}\rightarrow\alpha$. The continuity of $\Psi$ gives $\Psi(\alpha_{i}1_{P})\rightarrow \Psi(\alpha1_{P})$. So we obtain  that $\xi(\alpha_{i})\rightarrow\xi(\alpha)$. Then $\xi$ is continuous. Thereby, there exists a continuous function $\xi:\mathbb{R}\rightarrow\mathbb{R}$ such that
\begin{equation*}
\Psi(\alpha1_{P})=\xi(\alpha)\mathrm{K}(1_{P})+s\rho_{\frac{\pi}{2}}
\end{equation*}
for  $P\in\mathcal{P}^2$.
\qed

\begin{lem}\label{lem12}
Suppose that $n\geqslant3$. If the operator  $\Psi:L^{p}(\mu_{n})\rightarrow\langle\mathbb{M}^n,+\rangle$ is a continuous  $\mathrm{SL}(n)$ covariant valuation, then there exists  $\xi\in C(\mathbb R)$ such that
\begin{equation*}
\Psi(\alpha1_{P})=\xi(\alpha)\mathrm{K}(\alpha1_{P})
\end{equation*}
for  $P\in\mathcal{P}^n$.
\end{lem}

\emph{Proof.} For some $\alpha\in\mathbb{R}$, define $\mathrm{Y}_{\alpha}:\mathcal{P}^n\rightarrow\langle\mathbb{M}^n,+\rangle$ by setting
\begin{equation*}
\mathrm{Y}_{\alpha}(P)=\Psi(\alpha1_{P}).
\end{equation*}
Since $\Psi$ is a valuation in $L^{p}(\mu_{n})$, it follows for $P_{1},P_{2},P_{1}\cup P_{2},P_{1}\cap P_{2}\in\mathcal{P}^n$ that
\begin{align*}
\mathrm{Y}_{\alpha}(P_{1}\cup P_{2})+\mathrm{Y}_{\alpha}(P_{1}\cap P_{2}) &=\Psi(\alpha1_{P_{1}\cup P_{2}})+\Psi(\alpha1_{P_{1}\cap P_{2}})\\
&=\Psi(\alpha1_{P_{1}}\vee \alpha1_{P_{2}})+\Psi(\alpha1_{P_{1}}\wedge \alpha1_{P_{2}})\\
&=\mathrm{Y}_{\alpha}(P_{1})+\mathrm{Y}_{\alpha}(P_{2}).
\end{align*}
Therefore, $\mathrm{Y}_{\alpha}:\mathcal{P}^n\rightarrow \langle\mathbb{M}^n,+\rangle$ is a valuation.

\par  By Lemma \ref{lem2}, for $P\in\mathcal{P}^n$,   there is a sequence $P_{i}$ be a unions of finite closed cubes satisfying  $\lambda(P_{i}\bigtriangleup P)<1/i$. For $\alpha\in\mathbb{R}$, if $x\in P_{i}\cup P$, then there exists a constant $a>0$ such that $\vert x\vert\leq a$. We observe the following
\begin{align*}
\Vert\alpha1_{P_{i}}-\alpha1_{P}\Vert_{L^{p}(\mu_{n})}
 &=\vert\alpha\vert
\bigg(\int_{\mathbb{R}^n}\vert1_{P_{i}}(x)-1_{P}(x)\vert^{p}\vert x\vert^{2} dx\bigg)^{\frac{1}{p}}\\
 &\leq\vert\alpha\vert\bigg(\int_{\mathbb{R}^n}\vert1_{P_{i}}(x)-1_{P}(x)
\vert^{p}a^{2}dx\bigg)^{\frac{1}{p}}\\
 &=\vert \alpha\vert a^{\frac{2}{p}}\lambda(P_{i}\bigtriangleup P)^{\frac{1}{p}}.
\end{align*}
So if $\lambda(P_{i}\bigtriangleup P)\rightarrow 0$ as $i\rightarrow\infty$, then $\Vert\alpha1_{P_{i}}-\alpha1_{P}\Vert_{L^{p}(\mu_{n})}\rightarrow0$. The continuity of $\Psi$ gives $\Psi(\alpha1_{P_{i}})\rightarrow \Psi(\alpha1_{P})$. Hence, $\mathrm{Y}_{\alpha}(P_{i})\rightarrow \mathrm{Y}_{\alpha}(P)$. It leads that $\mathrm{Y}_{\alpha} $ is continuous.
Due to  $\mathrm{SL}(n)$ covariance of $\Psi$, for  $P\in\mathcal{P}^n$ and $\phi\in \mathrm{SL}(n)$
\begin{equation*}
\mathrm{Y}_{\alpha}(\phi P)=\Psi(\alpha1_{\phi P})=\Psi(\alpha1_{P}\circ\phi^{-1})=\phi \Psi(\alpha1_{p})\phi^{t},
\end{equation*}
then $\mathrm{Y}_{\alpha}$ is $\mathrm{SL}(n)$ covariant.

\par By Lemma \ref{lem4},  we have $\Psi(0)=\mathrm{\textbf{0}}$. Note that $\Vert\alpha1_{P}\Vert_{L^{p}(\mu_{n})}=0$ for all $P\in\mathcal{P}^n$ with dim $(P)< n$, therefore $\alpha1_{P}=0$ a.e. in $L^{p}(\mu_{n})$. The continuity of $\Psi$ implies that
\begin{equation}
\Psi(\alpha1_{P})=\Psi(0)=\mathrm{\textbf{0}}
\end{equation}
for every $P\in\mathcal{P}^n$ with  dim $(P)< n$. Thus
$\mathrm{Y}_{\alpha}(P)=\textbf{0}$
for every $P\in\mathcal{P}^n$ with  dim $(P)< n$, that is, $\mathrm{Y}_{\alpha}$ is a simple valuation.

 From Lemma \ref{lem9} and Lemma \ref{lem10},  there exist a functional $s_{\alpha}:\mathbb{R}\rightarrow\mathbb{R}$ and a constant $\alpha\in\mathbb{R}$ such that
\begin{equation*}
\Psi(\alpha1_{P})=s_{\alpha}\mathrm{K}(1_{P})
\end{equation*}
for every $P\in\mathcal{P}^n$. Define the function $\xi(\alpha)=s_{\alpha}$ for $\alpha\in\mathbb{R}$.
Next, we demonstrate the continuity of  $\xi$. Consider $\alpha\in\mathbb{R}$ and  a sequence $\{\alpha_{i}\}$  $\subseteq \mathbb{R}$ such that $\alpha_{i}\rightarrow\alpha$. We observe the following:
\begin{align*}
\Vert\alpha_{i}1_{P}-\alpha1_{P}\Vert_{L^{p}(\mu_{n})}=\vert\alpha_{i}
-\alpha\vert\left(\int_{\mathbb{R}^n}\vert1_{P}(x)\vert^{p}\vert x\vert^{2}dx\right)^{\frac{1}{p}}.
\end{align*}
Accordingly, $\Vert\alpha_{i}1_{P}-\alpha1_{P}\Vert_{L^{p}(\mu_{n})}\rightarrow0$ as $\alpha_{i}\rightarrow\alpha$. The continuity of $\Psi$ gives $\Psi(\alpha_{i}1_{P})\rightarrow \Psi(\alpha1_{P})$. So we have  $\xi(\alpha_{i})\rightarrow\xi(\alpha)$, that is, \,$\xi$ is continuous. Therefore, there exists  $\xi\in C(\mathbb R)$ such that
\begin{equation*}
\Psi(\alpha1_{P})=\xi(\alpha)\mathrm{K}(1_{P})
\end{equation*}
for $P\in\mathcal{P}^n$.
\qed
\par

The  lemma below is similar to Lemma 3.6 in Tsang's paper \cite{Tsang1}. Hence, the proof process is omitted in this context.

\begin{lem}\label{lem14}
Suppose that $\xi\in C(\mathbb R)$ and $\xi\not\equiv0$. If $\mathrm{K}(\xi\circ h)<\infty$ for every $h\in L^{p}(\mu_{n})$, then $\xi\in \widetilde{A}(\mathbb{R}).$
\end{lem}

\begin{lem}\label{lem13}
Suppose that $n\geq2$ and $\Psi_{1}, \Psi_{2}:L^{p}(\mu_{n})\rightarrow\langle\mathbb{M}^n,+\rangle$ are continuous, $\mathrm{SL}(n)$  covariant valuations. If
\begin{equation}
\Psi_{1}(\alpha1_{P})-\Psi_{1}(0)=\Psi_{2}(\alpha1_{P})-\Psi_{2}(0)
\end{equation}
for every $P\in\mathcal{P}^n$ and $\alpha\in\mathbb{R}^n$, then
\begin{equation*}
\Psi_{1}(h)-\Psi_{1}(0)=\Psi_{2}(h)-\Psi_{2}(0)
\end{equation*}
for every $h\in L^{p}(\mu_{n})$.
\end{lem}
\emph{Proof.} From the valuation property of  $\Psi_{1}, \Psi_2$, we have
\begin{equation}\label{3}
\Psi_{k}(h\vee0)+\Psi_{k}(h\wedge0)=\Psi_{k}(h)+\Psi_{k}(0),\ \  k=1, 2.
\end{equation}
Define $\widetilde{\Psi}_{k}:L^{p}(\mu_{n})\rightarrow\langle\mathbb{M}^n,+\rangle$
as
\begin{equation*}
\widetilde{\Psi}_{k}(h)=\Psi_{k}(h)-\Psi_{k}(0), \ \ k=1, 2.
\end{equation*}
Obviously, $\widetilde{\Psi}_{k}$ is a continuous SL(n) covariant valuation in $L^{p}(\mu_{n})$.
Due to Lemma \ref{lem3} and Lemma \ref{lem4}, we have $\widetilde{\Psi}_{k}(0)=\textbf{0}$. Therefore, (\ref{3}) can be rewritten as
\begin{equation*}
\widetilde{\Psi}_{k}(h\vee0)+\widetilde{\Psi}_{k}(h\wedge0)=\widetilde{\Psi}_{k}(h).
\end{equation*}
Thus, we just need to prove the following
\begin{equation}\label{4}
\widetilde{\Psi}_{1}(h)=\widetilde{\Psi}_{2}(h)
\end{equation}
for all $h\in L^{p}(\mu_{n})$ with $h\geq0$ and $h\leq0$.

\par Recall that simple functions of the form $\sum\limits_{i=1}^{m}\alpha_{i}1_{P_{i}}$  where $\alpha_{i}\in\mathbb{R}$ and $P_{i}\in\mathcal{P}^n$ with pairwise disjoint interiors, are dense in $L^{p}(\mu_{n})$. Because of
the continuity of $\widetilde{\Psi}_{1}$ and $\widetilde{\Psi}_{2}$, it is sufficient to prove Eq.
 (\ref{4}) for any simple function $h$ of the form $\sum\limits_{i=1}^{m}\alpha_{i}1_{P_{i}}$.  Here, $\alpha_{i}\in\mathbb{R}$ and $P\in\mathcal{P}^n$ have pairwise disjoint interiors.

\par In the case of  $h\geq0$, let
\begin{equation*}
h=\alpha_{1}1_{P_{1}}\vee\cdots\vee\alpha_{m}1_{P_{m}},\ \ \ 1 \leq i \leq m,
\end{equation*}
where the real numbers $\alpha_{i}\geq 0,$ and  $P_{i}\in\mathcal{P}$ have pairwise disjoint interiors.
Observe that  $\alpha_{k_{1}}1_{P_{k_{1}}}\wedge\cdots\wedge\alpha_{k_{j}}1_{P_{k_{j}}}=0$
a.e. in $L^{p}(\mu_{n})$ for all $2\leq j\leq m$ and $1\leq k_{1}<\cdots\leq k_{j}\leq m$.  Moreover, since $\widetilde{\Psi}_{k}$ is continuous in $L^{p}(\mu_{n})$, we have that
 \begin{equation}\label{5}
 \widetilde{\Psi}_{k}(\alpha_{k_{1}}1_{P_{k_{1}}}\wedge
 \cdots\wedge\alpha_{k_{j}}1_{P_{k_{j}}})=\textbf{0},\ \ k=1, 2.
 \end{equation}
From (\ref{5}) and the inclusion-exclusion principle,  we obtain that
\begin{align*}
\widetilde{\Psi}_{k}(h)
&=\widetilde{\Psi}_{k}(\alpha_{1}1_{P_{1}}\vee\cdots\vee\alpha_{m}1_{P_{m}})\\
&=\sum\limits_{j=1}\limits^{m}(-1)^{j-1}\sum\limits_{1\leq k_{1}<\cdots\leq k_{j}\leq m}\widetilde{\Psi}_{k}(\alpha_{k_{1}}1_{P_{k_{1}}}\wedge
\cdots\wedge\alpha_{k_{j}}1_{P_{k_{j}}})\\
&=\widetilde{\Psi}_{k}(\alpha_{1}1_{P_{1}})+\cdots+
 \widetilde{\Psi}_{k}(\alpha_{m}1_{P_{m}}).
\end{align*}
In the case of  $h\leq0,$ let
\begin{equation*}
h=\alpha_{1}1_{P_{1}}\wedge\cdots\wedge\alpha_{m}1_{P_{m}},\ \ \ 1\leq i\leq m,
\end{equation*}
where $\alpha_{i}\leq0,$ and  $P_{i}\in\mathcal{P}$ have pairwise disjoint interiors. Similarly,  we have
\begin{equation*}
\widetilde{\Psi}_{k}(h)=\widetilde{\Psi}_{k}(\alpha_{1}1_{P_{1}})+\cdots
 +\widetilde{\Psi}_{k}(\alpha_{m}1_{P_{m}}).
\end{equation*}
In both cases, we have
\begin{align*}
\widetilde{\Psi}_{1}(h)
&=\sum\limits_{i=1}\limits^{m}\widetilde{\Psi}_{1}(\alpha_{i}1_{P_{i}})\\
&=\sum\limits_{i=1}\limits^{m}(\Psi_{1}(\alpha_{i}1_{P_{i}})-\Psi_{1}(0))
=\sum\limits_{i=1}\limits^{m}(\Psi_{2}(\alpha_{i}1_{P_{i}})-\Psi_{2}(0))\\
&=\sum\limits_{i=1}\limits^{m}\widetilde{\Psi}_{2}(\alpha_{i}1_{P_{i}})
=\widetilde{\Psi}_{2}(h).
\end{align*}
Thus, the proof of this lemma is completed.
\qed

\textbf{Proof of Theorem \ref{theo3}}\
From $\phi\rho_{\frac{\pi}{2}}\phi^{t}=\rho_{\frac{\pi}{2}}$, then $\rho_{\frac{\pi}{2}}$ is $\mathrm{SL}(2)$ covariant for all $\phi\in \mathrm{SL}(2)$.
And by  Lemma \ref{lem7}, it shows that
\begin{equation*}
h\mapsto\mathrm{K}(\xi\circ h)+s\rho_{\frac{\pi}{2}}
\end{equation*}
is  a continuous $\mathrm{SL}(2)$ covariant valuation on $L^{p}(\mu_{n})$.
On the other hand, let $\Psi:L^{p}(\mu_{n})\rightarrow\mathbb{M}^2$ be a continuous  $\mathrm{SL}(2)$ covariant  valuation of $L^{p}(\mu_{n})$. From Lemma \ref{lem11}, Lemma \ref{lem14} and Lemma \ref{lem13}, it can be  proved the reverse statement.

\par\textbf{Proof of Theorem \ref{theo4}}\
Lemma \ref{lem7}
describes that
\begin{equation*}
h\mapsto\mathrm{K}(\xi\circ h)
\end{equation*}
determines a continuous $\mathrm{SL}(n)$ covariant valuation on $L^{p}(\mu_{n})$.
Additionally, Lemma \ref{lem12}, Lemma \ref{lem14} and Lemma \ref{lem13} ensure the reverse statement of Theorem \ref{theo4}.




\end{document}